\documentclass[11pt, english]{smfart}
\usepackage{graphpap}
\usepackage{amssymb}
\newtheorem{theoreme}{Theorem}
\renewcommand{\theequation}{\arabic{equation}}
\newtheorem{proposition}{Proposition}
\newtheorem{lemme}[proposition]{Lemma}
\newtheorem{corollaire}[proposition]{Corollary}
\newtheorem{definition}[proposition]{Definition}
\newtheorem{remarque}[proposition]{Remark}

\numberwithin{equation}{section}
\numberwithin{proposition}{section}

\def\Im{\textrm{Im}}
\def\Re{\textrm{Re}} 
\def\11{{\rm 1~\hspace{-1.4ex}l} }
\def\R{\mathbb R}
\def\C{\mathbb C}

\def\N{\mathbb N}

\def\O{\mathbb O}

\begin{document}
\title[Wave equation and invariant measures]
{Invariant measure for a three dimensional nonlinear wave equation }
\author{Nicolas Burq}
\address{D\'epartement de Math\'ematiques, Universit\'e Paris XI, 91 405  Orsay Cedex, 
France, and Institut Universitaire de France}
\email{nicolas.burq@math.u-psud.fr}
\author{Nikolay Tzvetkov}
\address{D\'epartement de Math\'ematiques, Universit\'e Lille I, 59 655 Villeneuve d'Ascq Cedex, France}
\email{nikolay.tzvetkov@math.univ-lille1.fr}
\begin{abstract} 
We study the long time behavior of the subcritical (subcubic) defocussing nonlinear wave equation on the 
three dimensional ball, for random data of low regularity. We prove that for a large set of radial initial data 
in $\cap_{s<1/2} H^s(B(0,1))$ the equation is (globally in time) well posed and we construct an invariant 
measure.
\end{abstract}
\begin{altabstract}
On \'etudie le comportement en grand temps de l'\'equation des ondes non lin\'eaire souscritique 
(sous cubique) d\'efocalisante dans la boule de dimension $3$, pour des donn\'ees initiales al\'eatoires. 
On d\'emontre que pour de nombreuses donn\'ees initiales radiales dans $\cap_{s<1/2} H^s( B(0,1))$ 
le probl\`eme est globalement bien pos\'e et on construit une mesure invariante par le flot
\end{altabstract}
\subjclass{ 35Q55, 35BXX, 37K05, 37L50, 81Q20 }
\keywords{nonlinear wave equation, eigenfunctions, dispersive equations, invariant measures}
\maketitle
%
\section{Introduction}
Consider the wave equation with  Dirichlet boundary condition
\begin{equation}\label{1}
(\partial_{t}^{2}-\mathbf{\Delta})w+|w|^{\alpha}w=0,\quad (w,\partial_t w)|_{t=0}=(f_1,f_2), 
\quad u\mid_{\R_t \times \partial \Theta} =0,\qquad \alpha>0
\end{equation}
with radial {\it real valued} initial data $(f_1,f_2)$ posed on the unit ball $\Theta$ of $\R^3$
defined by $\Theta\equiv (x\in\R^3:|x|<1)$. It is well-known that the functions
$$
e_{n}(r)\equiv \frac{\sqrt{2}\sin(\pi nr)}{r},\quad n=1,2,3,\cdots,
$$
where $r=|x|$ form an orthonormal bases of the Hilbert space of $L^2$ radial functions on $\Theta$.
Moreover $e_n$ are the radial eigenfunctions of the Laplace operator
$-\mathbf{\Delta}$ with 
Dirichlet boundary conditions, associated to eigenvalues 
$
z_n^2=(\pi n)^2.
$
We have the following statement.
\begin{theoreme}\label{thm1}
Suppose that $\alpha<2$. Let us fix a real number $\sigma$ such that 
\begin{equation}\label{restriction}
\max\big(0,\frac{\alpha-1}{\alpha}\big) <\sigma<\frac{1}{2}\,.
\end{equation}
Let
$((h_n(\omega),l_n(\omega))_{n=1}^{\infty}$ be a sequence of independent
standard real gaussians on a probability space $(\Omega,{\mathcal A},p)$.
Consider (\ref{1}) with initial data 
$$
f_1(r,\omega)=\sum_{n=1}^{\infty}\frac{ h_{n}(\omega)}{z_n}e_{n}(r),\quad
f_2(r,\omega)=\sum_{n=1}^{\infty} l_{n}(\omega)e_{n}(r)\,.
$$
Then almost surely in $\omega\in \Omega$ the problem (\ref{1}) has a unique global
solution 
$$u\in C(\R,H^\sigma_{rad}(\Theta)),$$
 $H^{\sigma}_{rad}(\Theta)$ being the Sobolev
space of index $\sigma$ of radial functions on $\Theta$. The uniqueness holds in the following sense :
for every $T>0$ there exists a Banach space $X_T$ continuously embedded in $C([-T,T],H^\sigma_{rad}(\Theta))$
such that the solution is unique in $X_T$. Furthermore we have the estimate
$
\|u(t)\|_{H^\sigma_{rad}(\Theta)} \leq C(f_1, f_2) \log(2+|t|)^{1/2}\,.
$ 
\end{theoreme}
\begin{remarque}
For every $\varepsilon>0$, the functions $(f_1,f_2)(r,\omega)$ belong a.s. to 
$H^{\frac{1}{2}-\varepsilon}_{rad}(\Theta)\times H^{-\frac{1}{2}-\varepsilon}_{rad}(\Theta)$ 
(see Lemma~\ref{lem.gauss}), but the probability of the event 
$$
\Big\{\omega\,:\,(f_1(r,\omega),f_2(r,\omega))\in 
H^{\frac{1}{2}}_{rad}(\Theta)\times H^{-\frac{1}{2}}_{rad}(\Theta)\Big\}
$$ 
is zero. Thus in the above statement, we obtain global solutions for data of 
low regularity. Such a regularity seems to be out of reach of the present deterministic methods. 
\end{remarque}
The map $\omega\mapsto (f_1(r,\omega),f_2(r,\omega))$ induces a Wiener measure on any Sobolev space 
of radial functions of regularity $<1/2$. It turns out that a measure absolutely continuous with respect 
to this Wiener measure is invariant under the global flow established in Theorem~\ref{thm1}
(see Theorem~\ref{thm2} below for a precise statement).
This measure invariance implies recurrence properties of the flow thanks to Poincar\'e's recurrence theorem.

The proof of Theorem~\ref{thm1} and Theorem~\ref{thm2} uses the Hamiltonian structure of the wave equation
(\ref{1}). We approximate (\ref{1}) by Hamiltonian ODE's and we obtain the solutions of (\ref{1}) as limits
of the solutions of these ODE's. We can ensure the passage to the limit thanks to a local well-posedness result
for (\ref{1}) for data of low regularity and the Liouville theorem for divergence free vector fields 
applied to the approximating ODE's.
In the local well-posedness argument we need to establish a Strichartz inequality for the wave equation,
posed on the disc and radial initial data. We hope that our elementary proof may be of independent interest
(see \cite{BLP} for Strichartz inequalities for the wave equation posed on a domain with boundary).
Our construction is inspired by the considerations in the works by Bourgain~\cite{Bo1,Bo2}, and the second author~\cite{Tz1,Tz2} in the context of the nonlinear 
Schr\"odinger equation (see also \cite{KS,Zh}
for works on invariant measures for the nonlinear 
Schr\"odinger equation). A difficulty we had to overcome is that for the wave
equation, in contrast with 
the nonlinear Schr\"odinger equation the $L^2$ norm is not conserved under the flow. This implies the
failure of the construction of the statistical ensemble of \cite{Tz1,Tz2} in the context of the wave equation. 
Here, we define a statistical ensemble which is invariant under the flow but for reasons quite different
from \cite{Tz1,Tz2} (in \cite{Tz1,Tz2} the $L^2$  norm conservation is important for the argument).
Let us also mention that the globalization argument presented here is simplified with respect to an analogous 
consideration in \cite{Tz1,Tz2}.

\section{Reduction of the problem}
For $\sigma\in\R$, we define $H^{\sigma}_{rad}(\Theta)$ as
$$
H^{\sigma}_{rad}(\Theta)\equiv\Big(\sum_{n=1}^{\infty}c_{n}e_{n},\,c_{n}\in\C\,:\,
\sum_{n=1}^{\infty} z_{n}^{2\sigma}|c_n|^2<\infty\Big)
$$
(the convergence of $\sum_{n=1}^{\infty}c_{n}e_{n}$ being apriori understood in ${\mathcal D}'(\Theta)$).
We can then equip $H^{\sigma}_{rad}(\Theta)$ with the natural complex Hilbert space structure.
In the case $\sigma=0$, we denote $H^{0}_{rad}(\Theta)$ by $L^2_{rad}(\Theta)$ and we have that 
the scalar product on $L^2_{rad}$ is defined by $\langle f,g\rangle=\int_{\Theta}f\bar{g}$.
Moreover $H^{\sigma}_{rad}(\Theta)$ and $H^{-\sigma}_{rad}(\Theta)$ are in a natural duality and we
we denote by $\langle\cdot,\cdot\rangle$ their pairing (in the case $\sigma=0$
we simply have the $L^2_{rad}$ scalar product).
For $\gamma\in\R$, we define the map $\sqrt{ - \mathbf{\Delta}}^{\gamma}$
acting as isometry from 
$H^{\sigma}_{rad}(\Theta)$ to
$H^{\sigma-\gamma}_{rad}(\Theta)$ by
$$
\sqrt{ - \mathbf{\Delta}}^{\gamma}\Big(\sum_{n=1}^{\infty}c_n e_n\Big)=
\sum_{n=1}^{\infty}z_{n}^{\gamma}c_n e_n\,.
$$
Clearly $\sqrt{ - \mathbf{\Delta}}^{\gamma_1+\gamma_2}=
\sqrt{ - \mathbf{\Delta}}^{\gamma_1}\circ\sqrt{ - \mathbf{\Delta}}^{\gamma_2}$
and $\sqrt{ - \mathbf{\Delta}}^{0}$ is the identity.
For $\gamma>0$ the map  $\sqrt{ - \mathbf{\Delta}}^{\gamma}$ is acting as 
``a differentiation'' while for $\gamma<0$, 
it is a smoothing operator. For $f\in L^2_{rad}(\Theta)$, we have
$\sqrt{ - \mathbf{\Delta}}^2(f)=-\mathbf{\Delta}(f)$, where $\mathbf{\Delta}$
is the Dirichlet self-adjoint realisation of the Laplacian. Moreover, for $f\in H^1_{rad}(\Theta)$, we have 
$$
\langle\mathbf{\Delta}(f),f\rangle=-\|\,\sqrt{ - \mathbf{\Delta}}(f)\|_{L^2(\Theta)}^{2}=
-\|\nabla f\|_{L^2(\Theta)}^{2}=
-\int_{0}^{1}|\partial_r f(r)|^{2}rdr,
$$
where $\nabla=(\partial_{x_1},\partial_{x_2},\partial_{x_3})$. 

Let us make some algebraic manipulations on (\ref{1}) allowing to write it as a first order in $t$
equation. Since $w_t$ is one derivative less regular than $w$ it is natural to set 
$v\equiv \sqrt{ - \mathbf{\Delta}}^{-1}(w_t)$ or equivalently 
$w_t=\sqrt{ - \mathbf{\Delta}}(v)$. If we set $u\equiv w+iv$ then we have that
$u$ solves the equation 
\begin{equation}\label{2}
(i\partial_{t}-\sqrt{ - \mathbf{\Delta}})u-\sqrt{ - \mathbf{\Delta}}^{-1}\big(|\Re(u)|^{\alpha}\Re(u)\big)
=0,\quad u|_{t=0}=u_{0},\quad u|_{\R\times\partial\Theta}=0,
\end{equation}
where $u_0=f_1+i\sqrt{ - \mathbf{\Delta}}^{-1}f_2$.
Therefore we have a correspondence between the solutions of (\ref{1}) with {\it real valued}
data $(f_1,f_2)\in H^{\sigma}_{rad}(\Theta)\times H^{\sigma-1}_{rad}(\Theta)$ and (\ref{2}) with data in 
$H^{\sigma}_{rad}(\Theta)$.
We are going to analyse (\ref{2}) with data $u_0$ in $H^{\sigma}_{rad}(\Theta)$. Then the real part of the
solutions of (\ref{2}) solve (\ref{1}) with $f_1=\Re(u_0)$ and $f_2=\sqrt{ - \mathbf{\Delta}}(\Im(u_0))$.

Let us next formally derive a conservation law for (\ref{2}). We will actually not use this conservation
law directly for (\ref{2}), we will only need it for a finite dimensional approximation of it. 
In order to highlight the algebraic computation, we make a formal computation in the context of (\ref{2}).
Let us write (\ref{2}) in the form
\begin{equation}\label{3}
iu_{t}+\sqrt{ - \mathbf{\Delta}}^{-1}\big(\mathbf{\Delta} u-|\Re(u)|^{\alpha}\Re(u)\big)=0\,.
\end{equation}
Since $\Im\langle \sqrt{ - \mathbf{\Delta}}^{-1}f,f\rangle=0$, taking the pairing of (\ref{3}) with 
$\mathbf{\Delta} u-|\Re(u)|^{\alpha}\Re(u)$ and taking the imaginary part gives that
$$
\frac{1}{2}\|\sqrt{ - \mathbf{\Delta}} (u)\|_{L^2(\Theta)}^{2}+
\frac{1}{\alpha+2}\|\Re(u)\|_{L^{\alpha+2}(\Theta)}^{\alpha+2}
$$
is conserved by the flow of (\ref{3}).

The free evolution associated to (\ref{3}) is given by the linear map $e^{-it\sqrt{ - \mathbf{\Delta}}}$ on 
$H^{\sigma}_{rad}(\Theta)$, $\sigma\in\R$, defined by
$$
e^{-it\sqrt{ - \mathbf{\Delta}}}\Big(\sum_{n=1}^{\infty}c_n e_n\Big)=\sum_{n=1}^{\infty}e^{-itz_{n}}c_n e_n\,.
$$
Observe that $e^{-it\sqrt{ - \mathbf{\Delta}}}$ acts as an isometry on $H^{\sigma}_{rad}(\Theta)$.
Let us also notice that thanks to the time oscillations for every $\sigma\in\R$ if 
$f\in H^{\sigma}_{rad}(\Theta)$ and $r_0\in(0,1]$ then
$e^{-it\sqrt{ - \mathbf{\Delta}}}(f)|_{\R\times\{r_0\}}$ is a well-defined distribution on $\R$.
In particular $e^{-it\sqrt{ - \mathbf{\Delta}}}(f)|_{\R\times\partial\Theta}=0$ 
and $e^{-it\sqrt{ - \mathbf{\Delta}}}(f)$ is the unique solution of
$(i\partial_{t}-\sqrt{ - \mathbf{\Delta}})u=0$ 
subject to the boundary condition
$u|_{\R\times\partial\Theta}=0$.
\section{Approximating ODE and associated gaussian measures}
Let us fix from now on a real number $\sigma$ satisfying (\ref{restriction}). 
Our analysis will be reduced to the study of
\begin{equation}\label{pak}
(i\partial_{t}-\sqrt{ - \mathbf{\Delta}})u-
\sqrt{ - \mathbf{\Delta}}^{-1}\big(|\Re(u)|^{\alpha}\Re(u)\big)=0,\quad
u|_{t=0}=u_{0},\quad u|_{\R\times\partial\Theta}=0,
\end{equation}
where the initial data $u_0$ belongs to $H^{\sigma}_{rad}(\Theta)$.
In order to prove Theorem~\ref{thm1}, we will need to study (\ref{pak}) with initial data given by
$$
u_{0}(r,\omega)=\sum_{n=1}^{\infty}\frac{ g_{n}(\omega)}{z_n}e_{n}(r),
$$
where $g_{n}(\omega)=h_{n}(\omega)+il_{n}(\omega)$ are independent normalized complex gaussian random variables.

For $N\geq 1$, we denote by $E_{N}$ the $N$ dimensional vector space on $\C$ spanned by $(e_n)_{n=1}^{N}$.
Let us denote by $S_{N}$ the projection on $E_{N}$ defined on every $H^{\sigma}_{rad}(\Theta)$ by
$$
S_{N}\Big(\sum_{n=1}^{\infty}c_{n}e_{n}\Big)\equiv\sum_{n=1}^{N}c_{n}e_{n}\,.
$$
We denote by $i_{N}$ the canonical isomorphism of vector spaces map from $\R^{2N}$ to $E_N$ defined by 
$$
i_{N}(((a_{n},b_n))_{n=1}^{N})\equiv\sum_{n=1}^{N}(a_{n}+ib_{n})e_{n}\,.
$$
The map $i_{N}$ equips $E_{N}$ with a canonical Borel sigma algebra and a canonical Lebesgue measure.

We shall approximate the solutions of (\ref{pak}) by the solutions of the ODE
\begin{equation}\label{pak_N}
(i\partial_{t}-\sqrt{ - \mathbf{\Delta}})u-
S_{N}\left(\sqrt{ - \mathbf{\Delta}}^{-1}\big(|\Re(u)|^{\alpha}\Re(u)\big)\right)=0,\,
u|_{t=0}=u_{0}\in E_{N},\, u|_{\R\times\partial\Theta}=0\,.
\end{equation}
Let us define the measure $\mu_{N}$ on $E_{N}$ 
as the image measure under the map from
$(\Omega,{\mathcal A},p)$ to $E_{N}$ (equipped with the Borel sigma algebra) defined by
\begin{equation}\label{zve}
\omega\longmapsto \sum_{n=1}^{N}\frac{h_{n}(\omega)+i l_{n}(\omega)}{z_n}e_{n}\,,
\end{equation}
where $h_{n}(\omega),l_{n}(\omega)$, $n=1,\cdots N$  is a sequence of
independent standard real gaussians ($h_n, l_n\in {\mathcal N}(0,1)$).
Observe that $\mu_N$ defines a probability measure on $E_N$.
We next define the measure $\rho_{N}$ as the image measure on $E_N$ by the map (\ref{zve})
of the measure
\begin{equation}\label{rho}
\exp\Big(-\frac{1}{(\alpha+2)}
\|\sum_{n=1}^{N}\frac{ h_{n}(\omega)}{z_n}e_{n}\|_{L^{\alpha+2}(\Theta)}^{\alpha+2}
\Big)dp(\omega).
\end{equation}
It turns out that $\rho_N$ is invariant under the flow of (\ref{pak_N}).
\begin{proposition}\label{liouville}
For every $u_0\in E_N$ the flow of (\ref{pak_N}) is defined globally in time.
Moreover the measure $\rho_N$ is invariant under this flow.
\end{proposition}
\begin{proof}
The local existence and uniqueness for the ODE (\ref{pak_N}) follows from the
Cauchy-Lipschitz theorem. Let us notice that the time existence given by the
Cauchy-Lipschitz theorem is very short (depending on $N$). We can however
extend globally in time the solutions of (\ref{pak_N}) thanks to the energy
conservation law associated to (\ref{pak_N}). Indeed if we multiply (\ref{pak_N})
by $\mathbf{\Delta} u-S_{N}(|\Re(u)|^{\alpha}\Re(u))$ (which is an element of $E_N$,
i.e. $C^{\infty}(\Theta)$ and vanishing on the boundary) and integrate over $\Theta$,
we get that the solutions of (\ref{pak_N}) satisfy
$$
\frac{d}{dt}
\Big[\frac{1}{2}\|\sqrt{ - \mathbf{\Delta}} (u)\|_{L^2(\Theta)}^{2}+
\frac{1}{\alpha+2}\|\Re(u)\|_{L^{\alpha+2}(\Theta)}^{\alpha+2}\Big]=0\,.
$$
Thus there exists a constant $C$ depending on $\sup_{1\leq n\leq N}|\langle
u_0,e_n\rangle|$ and $N$ but independent of $t$ such that as far as the
solution exists one has $\sup_{1\leq n\leq N}|\langle u(t),e_n\rangle|\leq C$.
Therefore the solutions of (\ref{pak_N}) are defined globally in time.
Let us now turn to the proof of the measure invariance.
Let us decompose the solution of (\ref{pak_N}) as
$$
u(t)=\sum_{n=1}^{N}\big(a_{n}(t)+ib_{n}(t)\big)e_{n},\quad a_{n}(t),b_n(t)\in\R\,.
$$
Then, if we set
$$
H(a_1,\dots,a_N,b_1,\dots,b_N)\equiv
\frac{1}{2}\sum_{n=1}^{N}z_n^2(a_n^2+b_n^2)+\frac{1}{\alpha+2}\int_{\Theta}
\big|\sum_{n=1}^{N}a_{n}e_{n}\big|^{\alpha+2}
$$
the problem (\ref{pak_N}) may be rewritten in the coordinates $a_n,b_n$ as
\begin{equation}\label{ab}
\dot{a}_{n}=z_{n}^{-1}\frac{\partial H}{\partial b_n},\quad 
\dot{b}_{n}=-z_{n}^{-1}\frac{\partial H}{\partial a_n},\quad n=1,\dots,N\,.
\end{equation}
Let us first observe that thanks to the structure of (\ref{ab}) the quantity $H(a_1,\dots,a_N,b_1,\dots,b_N)$
is conserved under the flow of (\ref{ab}). Let us also remark that
$$
\sum_{n=1}^{N}
\Big[
\frac{\partial}{\partial a_n}\big(z_{n}^{-1}\frac{\partial H}{\partial b_n}\big)+
\frac{\partial}{\partial b_n}\big(-z_{n}^{-1}\frac{\partial H}{\partial a_n}\big)
\Big]=0\,.
$$
Therefore we may apply Liouville's theorem for divergence free vector fields
to obtain that the measure 
$$
\prod_{n=1}^{N}da_n db_n
$$
is conserved by the flow of (\ref{ab}).
Since $ H(a_1,\dots,a_N,b_1,\dots,b_N)$ is conserved under the flow of
(\ref{ab}) we obtain that the measure
\begin{multline*}
\exp\big(-H(a_1,\dots,a_N,b_1,\dots,b_N)\big)\prod_{n=1}^{N}da_n
db_n
\\
=
\exp\Big(-\frac{1}{(\alpha+2)}\int_{\Theta}
\big|\sum_{n=1}^{N}a_{n}e_{n}\big|^{\alpha+2}\Big)
\prod_{n=1}^{N}
e^{-(z_n^2)(a_n^2/2)}da_n e^{-(z_n^2/)(b_n^2/2)}db_n
\end{multline*}
is also conserved by the flow of (\ref{ab}).
We therefore have that the measure
\begin{equation}\label{back}
\exp\Big(-\frac{1}{(\alpha+2)}\int_{\Theta}
\big|\sum_{n=1}^{N}a_{n}e_{n}\big|^{\alpha+2}\Big)
d\tilde{\mu}_{N},
\end{equation}
where 
$$
d\tilde{\mu}_{N}=(2\pi)^{-N}\big(\prod_{n=1}^{N}z_n^2\big)
\prod_{n=1}^{N}
e^{-(z_n^2)(a_n^2/2)}da_n e^{-(z_n^2)(b_n^2/2)}db_n
$$
is conserved by the flow of (\ref{ab}).
Observe that $d\tilde{\mu}_{N}$ is a probability measure on $\R^{2N}$.
The measure $d\tilde{\mu}_{N}$ is the distributions of the $\R^{2N}$ valued
random variable defined by
\begin{equation}\label{zve_bis}
\omega\longmapsto
\Big(
\frac{h_1(\omega)}{z_1},\frac{l_1(\omega)}{z_1},\cdots,
\frac{h_N(\omega)}{z_N},\frac{l_N(\omega)}{z_N}
\Big),
\end{equation}
where $h_n,l_n$, $n=1,\cdots N$ is again a system of independent standard real gaussians.
Moreover, 
the composition of the map (\ref{zve_bis}) and $i_{N}$ induces a probability measure
on $E_N$.
Next, coming back to (\ref{back}), we obtain that the image measure, induced on $\R^{2N}$ by the map
(\ref{zve_bis}), of the measure
\begin{equation*}
\exp\Big(-\frac{1}{(\alpha+2)}\int_{\Theta}
\big|\sum_{n=1}^{N}h_{n}(\omega)e_{n}\big|^{\alpha+2}\Big)
dp(\omega)
\end{equation*}
is invariant under the flow of (\ref{ab}).
This in turn implies that after applying $i_{N}$ we have a measure
on $E_N$ which is invariant under the flow of (\ref{pak_N}).
Coming back to (\ref{zve}), (\ref{rho}), we obtain that this measure is
precisely $\rho_N$. This completes the proof of Proposition~\ref{liouville}.
\end{proof}

Let us define the measure $\mu$ on $H^{\sigma}_{rad}(\Theta)$ (recall that
$\sigma$ is fixed and obeys (\ref{restriction})) as the image measure under the map from
$(\Omega,{\mathcal A},p)$ to $H^{\sigma}_{rad}(\Theta)$ equipped with the Borel sigma algebra, defined by
\begin{equation}\label{map}
\omega\longmapsto \sum_{n=1}^{\infty}\frac{ h_{n}(\omega)+i l_{n}(\omega)  }{z_n}e_{n}\,,
\end{equation}
where $((h_{n},l_{n}))_{n=1}^{\infty}$ is a sequence of independent standard real gaussians. 
Let us remark that the quantity 
$$
\sum_{n=1}^{\infty}\frac{h_{n}(\omega)+i l_{n}(\omega)  }{z_n}e_{n}
$$ 
is defined as the limit in $L^2(\Omega;H^{\sigma}_{rad}(\Theta))$ of the Cauchy sequence
$$
\sum_{n=1}^{N}\frac{ h_{n}(\omega)+i l_{n}(\omega)}{z_n}e_{n}
$$ 
and the mesurability of (\ref{map}) follows from the fact that 
the minimal sigma algebra containing the cylindrical sets of $H^{\sigma}_{rad}(\Theta)$
is the Borel sigma algebra.

Using \cite[Theorem~4]{AT}, we have that for $\alpha<4$ the quantity
$$
\|\sum_{n=1}^{\infty}\frac{ h_{n}(\omega)+i l_{n}(\omega)  }{z_n}e_{n}\|_{L^{\alpha+2}(\Theta)}
$$
is finite almost surely. 
Therefore, we can define a nontrivial measure $\rho$ on
$H^{\sigma}_{rad}(\Theta)$ as the image measure by the map (\ref{map}) of the measure
$$
\exp\Big(-\frac{1}{(\alpha+2)}
\|\sum_{n=1}^{\infty}\frac{ h_{n}(\omega)}{z_n}e_{n}\|_{L^{\alpha+2}(\Theta)}^{\alpha+2}
\Big)
dp(\omega).
$$
Observe that if a Borel set $A\subset H^{\sigma}_{rad}(\theta)$ is of full $\rho$ measure
then $A$ is also of full $\mu$ measure. Therefore, we need to solve (\ref{pak}) globally in time
for $u_0$ in a set of full $\rho$ measure.

We next turn to the limits of the measures $\rho_N$. As in \cite{Tz1,Tz2}, we can show
that if $U$ is an open set of $H^{s}_{rad}(\Theta)$, $s\in [\sigma,1/2[$ (and thus a Borel set of
$H^{\sigma}_{rad}(\Theta)$) then
\begin{equation}\label{open}
\rho(U)\leq\liminf_{N\rightarrow\infty}\rho_{N}(U\cap E_{N})\,.
\end{equation}
Moreover, if $F$ is a closed set of $H^{s}_{rad}(\Theta)$, $s\in [\sigma,1/2[$ then
\begin{equation}\label{close}
\rho(F)\geq\limsup_{N\rightarrow\infty}\rho_{N}(F\cap E_{N})\,.
\end{equation}
Let us remark that we have the following  standard gaussian estimate (see e.g. \cite{Tz1,Tz2}).
\begin{lemme}\label{lem.gauss}
Let $c$ be a positive constant satisfying $c<\pi/2$.
Denote by $B(0, \Lambda)_s$ the open ball of center $0$ and radius $\Lambda$ in $H^s_{rad}(\Theta)$. 
Then for every $s\in [\sigma,1/2[$, there exists $C_s>0$ 
such that for every $N, \Lambda$,
\begin{equation}\label{eq.gauss1}
\rho_{N}( B(0,\Lambda)_s^c\cap E_N)\leq \mu_N( B(0,\Lambda)^c\cap E_N)\leq C_s e^{-c\Lambda^2}\,.
\end{equation}
On the other hand, for every $s\geq \frac 1 2$ and every $N, \lambda$
\begin{equation}\label{eq.gauss2}
\rho_{N}( B(0,\Lambda)_s\cap E_N)\leq \mu_N( B(0,\Lambda)\cap E_N)\leq o(1)_{N\rightarrow + \infty}\,.
\end{equation}
As a consequence of~\eqref{open}, (\ref{close}) and~\eqref{eq.gauss1},~\eqref{eq.gauss2} we obtain
\begin{equation}\label{eq.gauss3}
\begin{gathered}
\rho( B(0,\Lambda)_s^c))\leq C_s e^{-c\Lambda^2}, \qquad s<\frac 1 2\\
\rho( B(0,\Lambda)_s))=0 \qquad s\geq\frac 1 2
\end{gathered}
\end{equation}
In particular for every $s<1/2$ the space $H^{s}_{rad}(\Theta)$ is of full $\rho$ measure but 
$\rho(H^{1/2}_{rad}(\Theta))=0$. 
\end{lemme}
\begin{proof}
The first inequality in~\eqref{eq.gauss1} is straightforward. The second is a
simple consequence of the Bienaym\'e-Tchebichev inequality.
More precisely, using the assumption on $c$ (recall that $z_n=\pi n$) 
and the Bienaym\'e-Tchebichev inequality yield
\begin{multline*}
e^{c\Lambda^2}
\mu_{N}( B(0,\Lambda)_s^c\cap E_N) \leq  
\int_{E_N} e^{c\|u\|_{H^s}^2} d\mu_N(u)
=
\prod_{n=1}^N\int_{\C} 
e^{cz_n^{2s}|c_n|^{2}-\frac{z_n^2|c_n|^2}{2}}\frac{z_n^2dc_c}{2\pi }
\\
 = 
\prod_{n=1}^N\int_{\C} e^{-\frac{|c_n|^2}{2} (1-\frac {2c}{z_n^{2-2s}})}
\frac{dc_n}{2\pi}=
\prod_{n=1}^N\frac{1}{1-\frac {2c}{z_n^{2-2s}}}\leq 
\prod_{n=1}^{\infty}\frac{1}{1-\frac {2c}{z_n^{2-2s}}}=C_{s},
\end{multline*}
where in the last inequality we used that $s<\frac 1 2$.

To prove~\eqref{eq.gauss2}, we again use  Bienaym\'e-Tchebichev inequality to write
$$
e^{-c\Lambda^2}
 \mu_N( B(0,\Lambda)\cap E_N)\
\leq 
\int_{E_N} e^{-c\|u\|_{H^s}^2} d\mu_N(u)
=
\prod_{n=1}^{N}\frac{1}{1+\frac {2c}{z_n^{2-2s}}}=
o_s(1)_{N\rightarrow + \infty},
$$
where in the last inequality we used that for $s\geq 1/2$, 
$$
\prod_{n=1}^{\infty}\frac{1}{1+\frac {2c}{z_n^{2-2s}}}=0\,.
$$
Finally, to prove the first part in~\eqref{eq.gauss3} we remark that 
$B(0,\Lambda)^c\subset \overline{B(0, \Lambda/2)}^c$ and apply~\eqref{open}
and~\eqref{eq.gauss1} and to prove the second part, we apply
directly~\eqref{open} and~\eqref{eq.gauss2}.
\end{proof}

\section{Strichartz estimates}
In this section we prove that the usual $3$-d Strichartz estimates are true
for our boundary value problem 
(and radial functions). Let us begin with a definition.
\begin{definition}
A couple of real numbers $(p,q), 2<p\leq + \infty$ is admissible if $\frac{1}{p}+\frac{1}{q}=\frac{1}{2}$.
For $T>0$, $0\leq s < 1$, we define the spaces 
$$ X^s_T= C^0 ([-T, T]; H^s_{rad}( \Theta)) \cap L^p((-T, T) ; L^q_{rad}( \Theta)), 
(p = \frac 2 s,q) \text{ admissible} 
$$
and it dual space
$$
 Y^s_T= L^1 ([-T, T]; H^{-s}_{rad}( \Theta)) + L^{p'}((-T, T) ; L^{q'}_{rad}( \Theta)), 
(p = \frac 2 s,q) \text{ admissible}
$$
equipped with their natural norms ($(p',q')$ being the conjugate couple of $(p,q)$).
\end{definition} 
\begin{proposition}\label{str}
Let $(p,q)$ be an admissible couple. 
Then there exists $C>0$ such that for every $T\in]0,1]$, 
every $f\in H^{\frac{2}{p}}_{rad}(\Theta)$ one has
\begin{equation}\label{eq.homog}
\|e^{-it\sqrt{- \mathbf{\Delta}}}(f)\|_{L^{p}([-T,T];L^{q}(\Theta))}\leq 
C\|f\|_{H^{\frac{2}{p}}_{rad}(\Theta)}\,.
\end{equation}
\end{proposition}
\begin{corollaire}\label{cor.stri}
For every $0< s <1 $, every admissible couple $(p,q)$, there exists $C>0$ such that
for every $T\in]0,1]$, every $f\in H^{\frac{2}{p}}_{rad}(\Theta)$ one has
\begin{equation}\label{eq.homogbis}
\|e^{-it\sqrt{- \mathbf{\Delta}}}(f)\|_{X^s_T}\leq 
C\|f\|_{H^{\frac{2}{p}}(\Theta)},\text{ if\, $\frac 1 p = \frac s 2$}
\end{equation}
\begin{equation}\label{eq.inhomog}
\|\int_0^t 
\sqrt{-\mathbf{\Delta}} ^{-1}
e^{-i (t-\tau)\sqrt{ - \mathbf{\Delta}}} (f)(\tau) d\tau\|_{X^s_T}\leq C \|f\|_{Y^{1-s}_T}
\end{equation}
\begin{equation}\label{eq.inhomogbis}
\|(1- S_N)\int_0^t\sqrt{-\mathbf{\Delta}} ^{-1} 
e^{-i (t-\tau)\sqrt{ - \mathbf{\Delta}}} (f)(\tau) d\tau\|_{X^s_T}
\leq C N^{s-s_1}\|f\|_{Y^{1- s_1}_T}, \text{ if $s<s_1<1 $}.
\end{equation}  
\end{corollaire}
\begin{proof}[Proof of Corollary~\ref{cor.stri}]
Inequality~\eqref{eq.homogbis} is obtained by using~\eqref{eq.homog} and the conservation of the $H^s$ norm. 
In order to prove~\eqref{eq.inhomog}, we set $K= e^{-it \sqrt {- \mathbf{\Delta}}}$. 
According to~\eqref{eq.homogbis}, $K$ is bounded from $H^s_{rad}$ to $X^s_T$.
Consequently $K^*$ is bounded from $Y^s$ to $H^{-s}_{rad}$.
Using the last property with $s$ replaced by $1-s$ (which remains in $]0,1[$ if $s\in]0,1[$)
and the fact that $\sqrt{- \mathbf{\Delta}} ^{-1} $ is bounded from $H^{s-1}_{rad}$ to $H^s_{rad}$, 
we obtain the following sequence of continuous mappings
\begin{equation}\label{reditza}
Y^{1-s}_{T}\stackrel{K^{\star}}{\longrightarrow}H^{s-1}_{rad}(\Theta)
\stackrel{ \sqrt{- \mathbf{\Delta}} ^{-1} }{\longrightarrow}H^{s}_{rad}(\Theta)
\stackrel{K}{\longrightarrow}X^{s}_{T}\,.
\end{equation}
On the other hand, it is easy to check that
$$ 
K\sqrt{- \mathbf{\Delta}} ^{-1}K^*(f) = 
\int_{-T}^T \sqrt{- \mathbf{\Delta}}^{-1}e^{-i (t-\tau) \sqrt{- \mathbf{\Delta}} } f(\tau) d\tau\,. 
$$
An application of Christ-Kisselev Lemma~\cite{ChKi} gives~\eqref{eq.inhomog} shows that the map
$$
f\longmapsto 
\int_{0}^t \sqrt{- \mathbf{\Delta}}^{-1}e^{-i (t-\tau) \sqrt{- \mathbf{\Delta}} } f(\tau) d\tau
$$
is bounded from $Y^{1-s}$ to $X^s$ which proves \eqref{eq.inhomog}.
Finally~\eqref{eq.inhomogbis} is obtained in a similar way replacing (\ref{reditza}) by
$$
Y^{1-s_1}_{T}\stackrel{K^{\star}}{\longrightarrow}H^{s_1-1}_{rad}(\Theta)
\stackrel{ \sqrt{- \mathbf{\Delta}} ^{-1} }{\longrightarrow}H^{s_1}_{rad}(\Theta)
\stackrel{ 1-S_N }{\longrightarrow}H^{s}_{rad}(\Theta)
\stackrel{K}{\longrightarrow}X^{s}_{T}
$$
and using that 
$(1- S_N) $ is bounded from $H^{s_1}$ to $H^s$ with norm
$\leq CN^{s-s_1}$.
\end{proof}
Let us come back to the proof of Proposition~\ref{str}. 
We first notice that it is enough to prove a similar result for the solutions of wave equations.
\begin{proposition}\label{stribis}
Let $(p,q)$ be an admissible couple. 
Then there exists $C>0$ such that for every $T\in]0,1]$, every 
$(u_0, u_1)\in H_{rad}^{\frac{2}{p}}(\Theta)\times H_{rad}^{\frac{2}{p}-1}(\Theta)$ and $u$ solution 
to the following wave equation
\begin{equation}\label{realwave} 
(\partial_t^2 - \mathbf{\Delta} )u =0, u\mid_{\partial \Theta}=0, u\mid_{t=0} = u_0, 
\partial_t u \mid_{t=0}= u_1,
\end{equation}
one has
$$
\|u\|_{L^{p}([-T,T];L^{q}(\Theta))}\leq C
\Big(\|u_0\|_{H_{rad}^{\frac{2}{p}}(\Theta)}+\|u_1\|_{H_{rad}^{\frac{2}{p}-1}(\Theta)}\Big)\,.
$$
\end{proposition}
Indeed, as 
$$ u= \cos(t\sqrt{- \mathbf{\Delta}})u_0 + 
\frac{\sin(t\sqrt{- \mathbf{\Delta}})} {\sqrt{- \mathbf{\Delta}}} u_1,
$$
Proposition~\ref{str} follows from Proposition~\ref{stribis} and the fact that 
$1/\sqrt{- \mathbf{\Delta}}$ is an isometry from
$H_{rad}^{\frac{2}{p}-1}(\Theta)$ to $ H_{rad}^{\frac{2}{p}}(\Theta)$. 
\begin{proof}[Proof of Proposition~\ref{stribis}]
Remark first that the bound given in Proposition~\ref{stribis} for $p= + \infty$ is trivial. 
As a consequence, it suffices to prove the bound for $p=2$, all other bounds following by interpolation. 
According to the finite speed of propagation for the solutions of wave equations, Proposition~\ref{stribis} 
is a local result which is known near any point in the interior of $\Theta$. Consequently is suffices 
to prove it near the boundary, replacing in the left hand side $L^q( \Theta)$ by 
$L^q( \{x ; |x| \in [1/2,1]\})$. But the conservation of the energy gives
$$ 
\int_{r=0}^1 |\partial_r u|^2(t,r) r^2 dr \leq 
\|u_0\|^2_{H_{rad}^{1}(\Theta)}+\|u_1\|^2_{L^2(\Theta)}
$$
and, the $r$ weight being irrelevant for $r \in [1/2, 1]$,  (one dimensional) Sobolev embedding gives
$$ \|u\|_{L^{\infty}([-1,1]_t\times [1/2,1])}\leq C\Big(\|u_0\|_{H^1_{rad}(\Theta)}+
\|u_1\|_{L^{2}_{rad}(\Theta)}\Big)\,
$$
which is stronger than the bounds given in Proposition~\ref{stribis} for $p=2$.
\end{proof}
\begin{remarque} The proof of the Strichartz estimate is much simplified by the radial assumption. 
However, the more general case of hyperbolic regime when one assume that the singularities of the wave 
are transversal to the boundary would still be true (but with a more technical proof involving 
parametrices constructed by reflections on the boundary).
\end{remarque}
\section{Local well-posedness}
If we set $S(t)=e^{-it\sqrt{- \mathbf{\Delta}}}$ then (\ref{2}) is reduced to the integral equation
\begin{equation}\label{Duhamel}
u(t)=S(t)u_0-
i\int_{0}^{t}S(t-\tau)\sqrt{- \mathbf{\Delta}}^{-1}\big(|\Re(u(\tau))|^{\alpha}\Re(u(\tau))\big)d\tau\,.
\end{equation} 
The next statement provides bounds on the right hand-side of (\ref{Duhamel}).
\begin{proposition}\label{calculus}
Let us fix $s$ such that
\begin{equation}\label{sigmas}
\max\big(0,\frac{\alpha -1 }{\alpha}\big)<s<\frac{1}{2}\,.
\end{equation}
Set $F(u)=|\Re (u)|^{\alpha}u$. Then  there exist $C>0$, $\delta >0$ such that for every $T\in]0,1]$, 
every $u,v\in X^s_T$, every $u_0\in
H^{s}_{rad}(\Theta)$, 
\begin{equation}\label{i}
\big\|S(t)u_0 \big\|_{X^s_T}\leq 
C\|u_0\|_{H^{s}_{rad}(\Theta)}\, ,
\end{equation}
\begin{equation}\label{ii}
\Big\|\int_{0}^{t}S(t-\tau)\sqrt{- \mathbf{\Delta}}^{-1}
F(u(\tau))d\tau\Big\|_{X^s_T}\leq CT^{\delta}\|u\|^{\alpha+1}_{X^s_T}
\end{equation}
\begin{equation}\label{ii'}
\Big\|(1- S_N)\int_{0}^{t}S(t-\tau)\sqrt{- \mathbf{\Delta}}^{-1}
F(u(\tau))d\tau\Big\|_{X^s_T}\leq CT^{\delta}N^{-\delta}\|u\|^{\alpha+1}_{X^s_T}\, ,
\end{equation}
\begin{equation}\label{iii}
\Big\|\int_{0}^{t}S(t-\tau)\sqrt{- \mathbf{\Delta}}^{-1}\Big(F(u(\tau))-F(v(\tau))\Big)d\tau\Big\|_{X^{s}_{T}}
\leq CT^{\delta}\Big(\|u\|^{\alpha}_{X^s_T}+
\|v\|^{\alpha}_{X^s_T}\Big)
\|u-v\|_{X^s_T}
\end{equation}
and 
\begin{multline}\label{iv}
\Big\|\int_{0}^{t}
S(t-\tau)\sqrt{- \mathbf{\Delta}}^{-1}
S_{N}\Big(F(u(\tau))-F(v(\tau))\Big)d\tau\Big\|_{X^{s}_{T}}
\\
\leq CT^{\delta}\Big(\|u\|^{\alpha}_{X^s_T}+\|v\|^{\alpha}_{X^s_T}\Big)\|u-v\|_{X^s_T}
\end{multline}
\end{proposition}
\begin{proof}
Estimate (\ref{i}) follows directly from Proposition~\ref{str}.
Let us next show (\ref{ii}).
According to Corollary~\ref{cor.stri} the left hand side of (\ref{ii}) is bounded by
\begin{equation}\label{gin}
C
\||\Re(u)|^{\alpha}\Re(u)\|_{L^{\widetilde{p}'}((-T, T); L^{\widetilde{q}'}(\Theta))}
\end{equation}
where $(\widetilde {p}= \frac{ 2} {1-s}, \widetilde {q})$ is an admissible couple. 
But 
\begin{equation*}
\frac {\alpha + 1 } p - \frac 1 {\widetilde {p}'} 
=\frac{(\alpha +1)s}{2}-\frac{1+s} {2} =  \frac {\alpha s } {2} -\frac 1 2
\end{equation*} 
and the conditions $s < \frac 1 2$ and $\alpha <2$ ensure that 
$$
\frac p  {\alpha + 1 }>\widetilde{p}' \,.
$$
On the other hand
\begin{equation*}
\frac {\alpha + 1 } q - \frac 1 {\widetilde{q}'} =
(\alpha + 1 ) \frac {(1-s)} 2- (1-\frac{ s} 2)= \frac{(\alpha -1)} 2  - \frac {\alpha s}{2} 
\end{equation*}
and the condition $\frac{ \alpha -1} \alpha <s$
ensures that
$$ \frac q {\alpha +1} > \widetilde{q}'\,.$$
As a consequence, since $\Theta$ is compact, 
applying H\"older inequality to~\eqref{gin}, we obtain (\ref{ii}).
To prove \eqref{ii'}, we simply remark that for $\sigma>s$ close enough to $s$, 
the admissible couple $$(\check {p}= \frac{ 2} {1-\sigma}, \check {q})$$ still satisfies
$$ \frac p  {\alpha + 1 }>\check{p}',\quad \frac q  {\alpha + 1 }>\check {p}'$$
and the same proof (using~\eqref{eq.inhomogbis} instead of ~\eqref{eq.inhomog}) gives~\eqref{ii'}.  
The proofs of (\ref{iii}) and (\ref{iv}) are very similar to that of (\ref{ii}) and will be omitted.
This completes the proof of Proposition~\ref{calculus}.
\end{proof}
As a consequence of Proposition~\ref{calculus}, we infer the following well-posedness results
for (\ref{pak}) and (\ref{pak_N}).
\begin{proposition}\label{lwp}
Let us fix $s$ satisfying (\ref{sigmas}). There exist $C>0$, $c\in]0,1]$, $\gamma>0$ such that for every 
$A>0$ if we set $T=c(1+A)^{-\gamma}$ then for every $u_0\in H^{s}_{rad}(\Theta)$ 
satisfying $\|u_0\|_{H^{s}_{rad}}\leq A$ there exists a unique
solution $u$ of (\ref{pak}) in $X^{s}_{T}$. 
Moreover $\|u\|_{X^s_T}\leq {C}\|u_0\|_{H^{s}_{rad}(\Theta)}$.
Finally if $u$ and $v$ are two solutions with data $u_0$, $v_0$ respectively,
satisfying 
$
\|u_0\|_{H^{s}_{rad}}\leq A $ and $\|v_0\|_{H^{s}_{rad}}\leq A$ then 
$
\|u-v\|_{X^s_T}\leq C\|u_0-v_0\|_{H^{s}_{rad}(\Theta)} .
$
\end{proposition}
Thanks to (\ref{iv}) we also have a well-posedness
in the context of (\ref{pak_N}) with bounds independent of $N$. 
\begin{proposition}\label{lwpbis}
Let us fix $s$ satisfying (\ref{sigmas}). There exist $C>0$, $c\in]0,1]$, $\gamma>0$ such that for every 
$A>0$ if we set $T=c(1+A)^{-\gamma}$ then for every $N\in \mathbb{N}$ and every $u_0\in H^{s}_{rad}(\Theta)\cap E_N$ 
satisfying $\|u_0\|_{H^{s}_{rad}}\leq A$ there exists a unique
solution $u=S_{N}(u)$ of 
(\ref{pak_N}) in $X^{s}_T$. Moreover 
$\|u\|_{X^s_T}\leq {C}\|u_0\|_{H^{s}_{rad}(\Theta)}$.
Finally if $u$ and $v$ are two solutions with data $u_0$, $v_0$ respectively,
satisfying $\|u_0\|_{H^{s}_{rad}}\leq A$ and $\|v_0\|_{H^{s}_{rad}}\leq A$ then 
$\|u-v\|_{X^s_T}\leq C\|u_0-v_0\|_{H^{s}_{rad}(\Theta)}$.
\end{proposition}
\section{Global existence for (\ref{pak}) on a set of full $\rho$ measure}
Recall that $\sigma$ is a fixed number satisfying (\ref{restriction}) and
the dependence on $\sigma$ of several numerical constants and sets appearing in the
sequel will not always be explicitly mentioned.
Let us denote by $\Phi_{N}(t):E_{N}\rightarrow E_{N}$, $t\in\R$ the flow of (\ref{pak_N}) defined in
Proposition~\ref{liouville}.
In the next proposition, we obtain a long time bound for the solutions of (\ref{pak_N}) in weak topologies.
Observe that bounds in terms of the $H^1$ norm of the data are trivial by the Hamiltonian conservation 
but insufficient for our purposes.
\begin{proposition}\label{longtime}
For every integer $i\geq 1$, every integer $N\geq 1$, there exists 
a $\rho_{N}$ measurable set $\Sigma_{N}^{i}\subset E_{N}$ such that
$
\rho_{N}(E_{N}\backslash \Sigma_{N}^{i})\leq 2^{-i}
$
and there exists a constant $C$ such that for every $i\in\N$, 
every $N\in\N$, every $u_0\in \Sigma_{N}^{i}$, every $t\in\R$, 
$$
\|\Phi_{N}(t)(u_0)\|_{H^{\sigma}_{rad}(\Theta)}\leq C(i+\log(1+|t|))^{\frac{1}{2}}\,.
$$
\end{proposition}
\begin{proof}
For $i,j$ integers $\geq 1$, we set
$$
B_{N}^{i,j}(D)\equiv
\big\{u\in E_{N}\,:\,\|u\|_{H^{\sigma}_{rad}(\Theta)}\leq D(i+j)^{\frac{1}{2}}\big\},
$$
where the number $D\gg 1$ (independent of $i,j,N$) will be fixed later. 
Thanks to Proposition~\ref{lwpbis}, there exist $c>0$, $C>0$, $\gamma>0$ only depending on $\sigma$ 
such that if we set $\tau \equiv c D^{-\gamma}(i+j)^{-\gamma/2}$ then for every $t\in[-\tau,\tau]$,
\begin{equation}\label{preser}
\Phi_{N}(t)\big(B_{N}^{i,j}(D)\big)\subset 
\Big\{u\in E_{N}\,:\,\|u\|_{H^{\sigma}_{rad}(\Theta)}\leq C\,D(i+j)^{\frac{1}{2}}\}\, .
\end{equation}
Next, we set
$$
\Sigma_{N}^{i,j}(D)\equiv
\bigcap_{k=-[2^{j}/\tau]}^{[2^{j}/\tau]}\Phi_{N}(-k\tau)(B_{N}^{i,j}(D))\, ,
$$
where $[2^{j}/\tau]$ stays for the integer part of $2^{j}/\tau$. 
Using the invariance of the measure $\rho_{N}$ by the flow $\Phi_{N}$ (Proposition~\ref{liouville}), 
we can write
\begin{equation*}
\rho_{N}(E_{N}\backslash\Sigma_{N}^{i,j}(D))
\leq 
(2[2^{j}/\tau]+1)\rho_{N}(E_N\backslash B_{N}^{i,j}(D))
\leq 
C2^{j}D^{\gamma}(i+j)^{\gamma/2}\rho_{N}(E_N\backslash B_{N}^{i,j}(D))\,.
\end{equation*}
Using Lemma~\ref{lem.gauss} we now deduce 
\begin{equation}\label{zvez}
\rho_{N}(E_{N}\backslash\Sigma_{N}^{i,j}(D))\leq
C2^{j}D^{\gamma}(i+j)^{\gamma/2}e^{-cD^2(i+j)}\leq 2^{-(i+j)},
\end{equation}
provided $D\gg 1$, independent of $i,j,N$.
Thanks to (\ref{preser}), we obtain that for
$u_0\in\Sigma_{N}^{i,j}(D)$, the solution of (\ref{pak_N}) with data $u_0$ satisfies
\begin{equation}\label{jjj1}
\|\Phi_{N}(t)(u_0)\|_{H^{\sigma}_{rad}(\Theta)}\leq CD(i+j)^{\frac{1}{2}},\quad |t|\leq 2^{j}\,.
\end{equation}
Indeed, for $|t|\leq 2^{j}$, we may find an integer $k\in [-[2^{j}/\tau],[2^{j}/\tau]]$ and 
$\tau_1\in [-\tau,\tau]$ so that $t=k\tau+\tau_1$ and thus 
$u(t)=\Phi_{N}(\tau_1)\big(\Phi_{N}(k\tau)(u_0)\big)$.
Since $u_0\in\Sigma_{N}^{i,j}(D)$ implies that $\Phi_{N}(k\tau)(u_0)\in 
B_{N}^{i,j}(D)$, we may apply (\ref{preser}) and arrive at (\ref{jjj1}).
Next, we set
$$
\Sigma_{N}^{i}=\bigcap_{j= 1}^{\infty}\Sigma_{N}^{i,j}(D)\,.
$$
Thanks to (\ref{zvez}),
$
\rho_{N}(E_{N}\backslash \Sigma_{N}^{i})\leq 2^{-i}\,.
$
In addition, using (\ref{jjj1}), we get that there exists $C$ such that for every $i$, every
$N$, every $u_0\in \Sigma_{N}^{i}$, every $t\in \R$,
$$
\|\Phi_{N}(t)(u_0)\|_{H^{\sigma}_{rad}(\Theta)}\leq C(i+\log(1+|t|))^{\frac{1}{2}}\,.
$$
Indeed for $t\in \R$ there exists $j\in\N$ such that $2^{j-1}\leq 1+|t|\leq 2^j$ and we apply 
(\ref{jjj1}) with this $j$.
This completes the proof of Proposition~\ref{longtime}.
\end{proof}
For  integers $i\geq 1$ and $N \geq 1$, we define the
cylindrical sets
$$
\tilde{\Sigma}_{N}^{i}\equiv
\big\{
u\in H^{\sigma}_{rad}(\Theta)\,:\, S_{N}(u)\in \Sigma_{N}^{i}
\big\}.
$$
Next, for an integer $i\geq 1$, we set
$$
\Sigma^{i}\equiv 
\big\{
u\in H^{\sigma}_{rad}(\Theta)\,:\,\exists\, N_k\rightarrow\infty, N_k\in\N,\,
\exists\, u_{N_k}\in \Sigma_{N_k}^{i},\,u_{N_k}\rightarrow u\,\, {\rm in}\, H^{\sigma}_{rad}(\Theta)
\big\}.
$$
Then the set $\Sigma^{i}$ is a closed set of $H^{\sigma}_{rad}(\Theta)$.
Observe that we have the inclusion
$$
\limsup_{N\rightarrow\infty}\tilde{\Sigma}_{N}^{i}
\equiv\bigcap_{N=
  1}^{\infty}\bigcup_{N_1=N}^{\infty}\tilde{\Sigma}_{N_1}^{i}
\subset \Sigma^{i}.
$$
Therefore 
\begin{equation}\label{kr1}
\rho(\Sigma^{i})   \geq   \rho(\limsup_{N\rightarrow\infty}\tilde{\Sigma}_{N}^{i})\,.
\end{equation}
Using Fatou's lemma, we get
\begin{equation}\label{kr2}
\rho(\limsup_{N\rightarrow\infty}\tilde{\Sigma}_{N}^{i})
 \geq 
\limsup_{N\rightarrow \infty}\rho(\tilde{\Sigma}_{N}^{i})\,.
\end{equation}
Next, using  Proposition~\ref{longtime} and (\ref{open}) and (\ref{close}), we obtain
\begin{equation}\label{kr4}
\limsup_{N\rightarrow \infty}\rho(\tilde{\Sigma}_{N}^{i})
=
\limsup_{N\rightarrow \infty}\rho_{N}(\Sigma_{N}^{i})
\geq
\limsup_{N\rightarrow \infty}\big(\rho_{N}(E_{N})-2^{-i}\big)
=
\rho\big(H^{\sigma}_{rad}(\Theta)\big)-2^{-i}.
\end{equation}
Collecting (\ref{kr1}), (\ref{kr2}) and (\ref{kr4}), we arrive at
$$
\rho(\Sigma^{i})  \geq \rho\big(H^{\sigma}_{rad}(\Theta)\big)-2^{-i}.
$$
Now, we set
$$
\Sigma\equiv\bigcup_{i\geq 1}\Sigma^{i}\,.
$$
Thus $\Sigma$ is of full $\rho$ measure.
It turns out that one has global existence for $u_0\in \Sigma$. 
\begin{proposition}\label{global_existence}
Let us fix $i\in\N$.
Then for every $u_0\in \Sigma^{i}$, the local solution $u$ 
of (\ref{pak}) given by Proposition~\ref{lwp} is globally defined.
In addition there exists $C>0$ such that for every $u_0\in \Sigma^{i}$,
\begin{equation}\label{growth}
\|u(t)\|_{H^{\sigma}_{rad}(\Theta)}\leq C(i+\log(1+|t|))^{\frac{1}{2}}\,.
\end{equation}
Moreover, if $(u_{0,k})_{k\in\N}$, $u_{0,k}\in \Sigma^{i}_{N_k}$,
$N_k\rightarrow\infty$ converges to $u_0$ as $k\rightarrow\infty$ in $H^{\sigma}_{rad}(\Theta)$ then  
for every $t\in\R$,
\begin{equation}\label{limit}
\lim_{k\rightarrow\infty}\|u(t)-\Phi_{N_k}(t)(u_{0,k})\|_{H^{\sigma}_{rad}(\Theta)}=0\,.
\end{equation}
\end{proposition}
\begin{proof}
Let $u_0\in \Sigma^{i}$ and $u_{0,k}\in \Sigma^{i}_{N_k}$,
$N_k\rightarrow\infty$ a sequence
tending to $u_0$ in $H^{\sigma}_{rad}(\Theta)$. Let us fix $T>0$. Our aim so to extend the solution 
of (\ref{pak}) given by Proposition~\ref{lwp} to the interval $[-T,T]$.
Using Proposition~\ref{longtime}, we have that there exists a constant $C$ such that for every $k\in\N$,
every $t\in\R$,
\begin{equation}\label{ant}
\|\Phi_{N_k}(t)(u_{0,k})\|_{H^{\sigma}_{rad}(\Theta)}\leq C(i+\log(1+|t|))^{\frac{1}{2}}\,.
\end{equation}
Therefore, if we set $u_{N_k}(t)\equiv  \Phi_{N_k}(t)(u_{0,k})$ and 
$\Lambda\equiv C(i+\log(1+T))^{\frac{1}{2}}$, we have the bound
\begin{equation}\label{david1}
\|u_{N_k}(t)\|_{H^{\sigma}_{rad}(\Theta)}\leq \Lambda,\quad \forall\,|t|\leq T,\quad \forall\, k\in\N.
\end{equation}
In particular $\|u_0\|_{H^{\sigma}_{rad}}\leq\Lambda$ 
(apply (\ref{david1}) with $t=0$ and let $k\rightarrow\infty$).
Let $\tau>0$ be the local existence time for (\ref{pak}), provided by
Proposition~\ref{lwp} for $A=\Lambda+1$. 
Recall that we can assume $\tau=c(1+\Lambda)^{-\gamma}$
for some $c>0$, $\gamma>0$ depending only on the choice of $\sigma$. 
We can assume that $T>\tau$.
Denote by $u(t)$ the solution of (\ref{pak}) with data $u_0$ on the time interval $[-\tau,\tau]$. Then
$v_{N_k}\equiv u-u_{N_k}$ solves the equation
\begin{equation}\label{eqnv}
(i\partial_{t}-\sqrt{- \mathbf{\Delta}}) v_{N_k} =
\sqrt{- \mathbf{\Delta}}^{-1}\Big(F(u)-S_{N_k}(F(u_{N_k}))\Big), \quad v_{N_k}|_{t=0}=u_0-u_{0,k} \, ,
\end{equation}
where $F(u)=|\Re(u)|^{\alpha}\Re(u)$. Next, we write
$$
F(u)-S_{N_k}(F(u_{N_k}))=S_{N_k}\big(F(u)-F(u_{N_k})\big)+(1-S_{N_k})F(u).
$$Therefore
\begin{multline*}
v_{N_k}(t)=S(t)(u_0-u_{0,k})
\\
-i\int_{0}^{t}S(t-\tau)\sqrt{- \mathbf{\Delta}}^{-1}S_{N_k}\big(F(u(\tau))-F(u_{N_k}(\tau))\big)d\tau
\\
-i\int_{0}^{t}S(t-\tau)\sqrt{- \mathbf{\Delta}}^{-1}(1-S_{N_k})F(u(\tau))d\tau\,.
\end{multline*}
Using Proposition~\ref{calculus}, we obtain that there exist $C>0$ and $\theta, \delta>0$ (depending only on $\sigma$) such that one has the bound
\begin{equation*}
\|(1- S_N)\int_{0}^{t}S(t-\tau)\sqrt{-\mathbf{\Delta}}^{-1}
F(u(\tau))d\tau\|_{X^{\sigma}_{\tau}}\leq C \tau^{\theta}N^{- \delta}\|u\|_{X^{\sigma}_{\tau}}
\big(1+\|u\|_{X^{\sigma}_{\tau}}^{\alpha}\big). 
\end{equation*}
Another use of Proposition~\ref{calculus} yields
\begin{equation*}
\|v_{N_k}\|_{X^{\sigma}_{\tau}}\leq 
C\Bigl(\|u_0-u_{0,k}\|_{H^{\sigma}_{rad}(\Theta)}
 + \tau^{\theta}\|v_{N_k}\|_{X^{\sigma}_{\tau}}
\big(1+\|u\|_{X^{\sigma}_{\tau}}^{\alpha}+
\|u_{N_k}\|_{X^{\sigma}_{\tau}}^{\alpha}
\big)\Bigr) + o(1)_{k\rightarrow + \infty}\,.
\end{equation*}
A use of Proposition~\ref{lwp} and Proposition~\ref{lwpbis} yields
\begin{multline*}
\|v_{N_k}\|_{X^{\sigma}_{\tau}} \\
\leq
C\|u_0-u_{0,k}\|_{H^{\sigma}_{rad}(\Theta)}+C\tau^{\theta}\|v_{N_k}\|_{X^{\sigma}_{\tau}}
\Big(
1+C\|u_{0}\|_{H^{\sigma}_{rad}(\Theta)}^{\alpha}+C\|u_{0,k}\|_{H^{\sigma}_{rad}(\Theta)}^{\alpha}\Big)
+
o(1)_{k\rightarrow + \infty}
\\
\leq 
C\|u_0-u_{0,k}\|_{H^{\sigma}_{rad}(\Theta)}
+
C\tau^{\theta}(1+\Lambda)^{\alpha}\|v_{N_k}\|_{X^{\sigma}_{\tau}} + o(1)_{k\rightarrow + \infty}\,.
\end{multline*}
Recall that $\tau=c(1+\Lambda)^{-\gamma}$, where $c>0$ and $\gamma>0$ are depending only on 
$\sigma$. In the last estimate the constants $C$ and $\theta$ also depend only
on $\sigma$. 
Therefore, if we assume that $\gamma>\alpha/\theta$ then the restriction on $\gamma$
remains to depend only on $\sigma$. Similarly, if we assume that $c$ is so small that
$
C\tau^{\theta}(1+\Lambda)^{\alpha}\leq
Cc^{\theta}(1+\Lambda)^{-\gamma\theta}(1+\Lambda)^{\alpha}
\leq Cc^{\theta}
<1/2
$
then the smallness restriction on $c$ remains to depend only on $\sigma$. 
Therefore, we have that after possibly slightly modifying the values of $c$ and $\gamma$ 
(keeping $c$ and $\gamma$ only depending on $\sigma$  and independent of $N_k$)
in the definition of $\tau$ that
\begin{eqnarray}\label{eq.est}
\|v_{N_k}\|_{X^{\sigma}_{\tau}}
\leq 
C\|v_{N_k}(0)\|_{H^{\sigma}_{rad}(\Theta)}+o(1)_{k\rightarrow + \infty}.
\end{eqnarray}
and passing to the limit in~\eqref{eq.est}, we obtain
\begin{equation*}
\lim_{k\rightarrow + \infty}
\|v_{N_k}\|_{L^\infty( [0,\tau];H^{\sigma}_{rad}(\Theta))}=0,
\end{equation*}
where $ \tau= c(1+\Lambda)^{-\gamma}$
and the constants $c$ and $\gamma$ depend only on $\sigma$.
Thus, via a use of
the triangle inequality,
\begin{equation}\label{david2}
\|u(t)\|_{H^{\sigma}_{rad}(\Theta)}\leq
\limsup_{k\rightarrow + \infty}\|u_{N_k}(t)\|_{H^{\sigma}_{rad}(\Theta)}\leq\Lambda,\quad |t|\leq\tau.
\end{equation}
In particular, we deduce
$$ \|u(\tau)\|_{H^{\sigma}_{rad}(\Theta)}\leq\Lambda
$$ and we can repeat the argument for obtaining (\ref{david2}) on 
$(\tau,2\tau)$, $(2\tau, 3\tau)$, ...$([\frac T \tau] \tau,([\frac T \tau]+1) \tau)$
(and similarly for negative times), giving (\ref{growth}) and (\ref{limit}).
This completes the proof of Proposition~\ref{global_existence}.
\end{proof}
Therefore we solved globally in time, with a suitable uniqueness, the problem
(\ref{pak}) on a set of full $\rho$ measure.
This completes the proof of Theorem~\ref{thm1}.
\section{Invariance of the measure $\rho$ }
Set
\begin{multline*}
\Sigma^{i}(M)\equiv 
\big\{
u\in H^{\sigma}_{rad}(\Theta)\,:
\,\exists\, \tau_k\in\R, |\tau_k|\leq M,
\\
\,\exists\, N_k\rightarrow\infty, N_k\in\N,\,
\exists\, u_{N_k}\in \Sigma_{N_k}^{i},\,
\Phi_{N_k}(\tau_k)u_{N_k}\rightarrow u\,\, {\rm in}\, H^{\sigma}_{rad}(\Theta)
\big\}.
\end{multline*}
The set $\Sigma^{i}(M)$ is a closed set of $H^{\sigma}_{rad}(\Theta)$.
Observe that $\Sigma_{\sigma}^{i}(0)$ is the set $\Sigma^{i}$ used in the 
proof of Theorem~\ref{thm1}. Next, we set
$$
{\bf \Sigma}^{i}=\bigcup_{M=1}^{\infty}\Sigma^{i}(M).
$$
The set ${\bf \Sigma}^{i}$ is $\rho$ measurable and
\begin{equation}\label{chti}
\rho({\bf \Sigma}^{i})
\geq
\rho(\Sigma^{i}(0))  \geq \rho\big(H^{\sigma}_{rad}(\Theta)\big)-2^{-i}.
\end{equation}
Proposition~\ref{global_existence} naturally extends to the set ${\bf \Sigma}^{i}$
\begin{proposition}\label{global_existence_bis}
Let us fix $M,i\in\N$.
Then for every $u_0\in {\bf \Sigma}^{i}(M)$, the local solution $u$ 
of (\ref{pak}) given by Proposition~\ref{lwp} is globally defined.
In addition there exists $C>0$ such that for every $u_0\in {\bf \Sigma}^{i}(M)$,
\begin{equation}\label{growth_bis}
\|u(t)\|_{H^{\sigma}(\Theta)}\leq C(i+\log(1+|M|+|t|))^{\frac{1}{2}}\,.
\end{equation}
Moreover, if $(u_{0,k})_{k\in\N}$, $u_{0,k}\in \Sigma^{i}_{N_k}$,
$N_k\rightarrow\infty$, $|\tau_k|\leq M$ are such that $\Phi_{N_k}(\tau_k)(u_{0,k})$
converges to $u_0$ as $k\rightarrow\infty$ in $H^{\sigma}_{rad}(\Theta)$ then  
\begin{equation}\label{limit_bis}
\lim_{k\rightarrow\infty}\|u(t)-\Phi_{N_k}(t+\tau_k)(u_{0,k})\|_{H^{\sigma}_{rad}(\Theta)}=0\,.
\end{equation}
\end{proposition}
Next, we set
$$
{\bf \Sigma}\equiv \bigcup_{i=1}^{\infty}{\bf \Sigma}^{i}\,.
$$
Then, using (\ref{chti}) we obtain that the set ${\bf \Sigma}$ is of full $\rho$ measure.
Thanks to Proposition~\ref{global_existence_bis}, we can establish a well-defined dynamics of (\ref{pak})
for data in ${\bf \Sigma}$. Let us denote by $\Phi$ the flow map of (\ref{pak}) for data in
${\bf \Sigma}$. We have the following corollary of Proposition~\ref{global_existence_bis} and
Proposition~\ref{lwp}.
\begin{proposition}\label{invariance}
For every $t\in\R$, $\Phi(t)({\bf \Sigma} )={\bf \Sigma}$.
In addition $\Phi(t)$ is continuous with respect to the induced by $H^\sigma_{rad}(\Theta)$ 
to ${\bf \Sigma}$ topology (in particular $\rho$ measurable).
\end{proposition}
We now state the measure invariance result.
\begin{theoreme}\label{thm2}
For every $A\subset {\bf \Sigma} $, a $\rho$ measurable set, for every $t\in\R$,
$\rho(A)=\rho(\Phi(t)(A))$.
\end{theoreme}
\begin{proof}
We first perform several reductions allowing to consider only sets $A$ of a special type and only short 
times $t$. Thanks to the invariance of ${\bf \Sigma}$ under $\Phi(t)$ and the time reversibility
of  $\Phi(t)$, we obtain that it suffices to prove that for every $\rho$ measurable set 
$A\subset {\bf \Sigma}$, every $t\geq 0$
one has $\rho(A)\leq \rho(\Phi(t)(A))$. Next, we observe that every $\rho$ measurable set $A$ may be 
approximated from the interior by closed sets of $H^\sigma_{rad}(\Theta)$, 
i.e. there exists a sequence of closed
sets $F_n\subset A$ such that $\rho(A)=\lim_{n}\rho(F_n)$. Indeed, this approximation property is equivalent to
a similar approximation property from the exterior by open sets which may be achieved by considering 
$\varepsilon$ open neighborhoods of the set a passing to the limit $\varepsilon\rightarrow 0$ 
via the Lebesgue dominated convergence theorem.
Therefore, we deduce that it suffices to prove that for every closed set $F\subset {\bf \Sigma}$
of $H^\sigma_{rad}(\Theta)$
one has $\rho(F)\leq \rho(\Phi(t)(F))$, $t\geq 0$. 
Indeed, if we have the last inequality then for an arbitrary measurable set $A\subset {\bf \Sigma}$, 
we may write
$$
\rho(A)=\lim_{n\rightarrow\infty}\rho(F_n)\leq \limsup_{n\rightarrow\infty}\rho(\Phi(t)(F_n))\leq
\rho(\Phi(t)(A)),
$$
where $F_n\subset A$ is the corresponding approximating sequence of closed sets.
Let $F$ be a closed set of $H^\sigma_{rad}(\Theta)$. Let us consider the set $K_{n}\subset A$ defined
as 
$$
K_{n}\equiv \{u\in F\,:\, \|u\|_{H^s_{rad}(\Theta)}\leq n\},
$$
where $\sigma<s<1/2$. Then $K_{n}$ is a compact set of $H^{\sigma}_{rad}(\Theta)$ and thanks to
Lemma~\ref{lem.gauss} one has $\rho(F)=\lim_{n}\rho(K_n)$.
Therefore, in order to prove Theorem~\ref{thm2}, it suffices to prove that for every set 
$K\subset {\bf \Sigma}$ which is 
a compact of $H^{\sigma}_{rad}(\Theta)$ one has $\rho(K)\leq \rho(\Phi(t)(K))$, $t\geq 0$.

Let us now fix a compact $K\subset {\bf \Sigma}$ of $H^{\sigma}_{rad}(\Theta)$ and $t\geq 0$.
Let us observe that there exists $R>0$ such that
$$
\{\Phi(\tau)(K),\,\, 0\leq \tau\leq t\}\subset
\{u\in H^{\sigma}_{rad}(\Theta)\,:\, \|u\|_{H^\sigma_{rad}(\Theta)}\leq R\}\equiv B_{R}\,.
$$
We next state a proposition which allows to compare 
$\Phi$ and $\Phi_N$ for data in compacts contained in $B_{R}$.
\begin{lemme}\label{nedelia}
There exist two constants $c>0$ and $\gamma>0$ (depending only on $\sigma$) such that the following holds true.
For every compact $K\subset B_{R}$,
every $\varepsilon>0$ there exists $N_0\geq 1$ such that for every $N\geq N_0$, every $u_0\in K$,
every $\tau\in [0,c(1+R)^{-\gamma}]$,
$$
\|\Phi(\tau)(u_0)-\Phi_{N}(\tau)(S_{N}(u_0))\|_{H^{\sigma}_{rad}(\Theta)}<\varepsilon\,.
$$
\end{lemme}
\begin{proof}
The argument is very similar to Proposition~\ref{global_existence}, the only additional point is the 
uniformness with respect to the compact $K$, we will use below.
For $u_0\in K$, 
we denote by $u$ the solution of (\ref{pak}) with data $u_0$
and by $u_{N}$ the solution of (\ref{pak_N}) with data $S_{N}(u_0)$, defined on $[0,\tau]$, where
thanks to Proposition~\ref{lwp} and Proposition~\ref{lwpbis}.
$\tau=c_0(1+R)^{-\gamma_0}$ with $c_0>0$, $\gamma_0>0$ depending only on $\sigma$.
Next, we set $v_{N}\equiv u-u_N$. Then $v_{N}$ solves
\begin{equation}\label{eqnvpak}
(i\partial_{t}-\sqrt{- \mathbf{\Delta}}) v_{N} =
\sqrt{- \mathbf{\Delta}}^{-1}\Big(F(u)-S_{N}(F(u_{N}))\Big), \quad v_{N}(0)=(1-S_{N})u_0\,.
\end{equation}
where $F(u)=|\Re(u)|^{\alpha}\Re(u)$.
By writing
$$
F(u)-S_{N}(F(u_{N}))=S_{N}\big(F(u)-F(u_{N})\big)+(1-S_{N})F(u)
$$
and using Proposition~\ref{calculus}, we obtain that there exist $C>0$ and $\theta>0$ 
depending only on $\sigma$ such that 
\begin{eqnarray*}
\|v_{N}\|_{X^{\sigma}_{\tau}}
\leq C\|(1-S_N)u_0\|_{H^{\sigma}_{rad}(\Theta)} + C\tau^{\theta}\|v_N\|_{X^{\sigma}_{\tau}}
\big(1+\|u\|_{X^{\sigma}_{\tau}}^{\alpha}+\|u_{N}\|_{X^{\sigma}_{\tau}}^{\alpha}\big)+
o(1)_{N\rightarrow + \infty},
\end{eqnarray*}
where $o(1)_{N\rightarrow + \infty}$ is a quantity which tends to zero as $N\rightarrow + \infty$,
{\it uniformly} with respect to $u_0\in K$.
Using Proposition~\ref{lwp} and Proposition~\ref{lwpbis}, we get
\begin{eqnarray*}
\|v_{N}\|_{X^{\sigma}_{\tau}}
\leq C\|(1-S_N)u_0\|_{H^{\sigma}_{rad}(\Theta)} + C\tau^{\theta}\|v_N\|_{X^{\sigma}_{\tau}}
\big(1+\|u_0\|_{H^{\sigma}_{rad}(\Theta)}^{\alpha}\big)+o(1)_{N\rightarrow + \infty}.
\end{eqnarray*}
Coming back to the definition of $\tau$
we can choose $c_0$ small enough and $\gamma_0$ large enough,
but keeping their dependence only on $\sigma$, to infer that
$$
\|v_N\|_{X^{\sigma}_{\tau}}\leq C\|(1-S_N)u_0\|_{H^{\sigma}_{rad}(\Theta)}\,.
$$
The space $X^{\sigma}_{\tau_0}$ is continuously 
embedded in $C([0,\tau];H^s_{rad}(\Theta))$ and thus there exists $C$ depending only on 
$\sigma$ such that
\begin{equation*}
\|v_N(t)\|_{H^{s}(\Theta)}\leq C\|(1-S_N)u_0\|_{H^{\sigma}_{rad}(\Theta)},\quad t\in [0,\tau].
\end{equation*}
Since $K$ is a compact of $H^{\sigma}_{rad}(\Theta)$, we have 
$$
\forall\,\varepsilon>0,\quad \exists N_0\geq 1\,:\, \forall\, N\geq
N_0,\,\forall\, u_0\in K,\,
\|(1-S_N)u_0\|_{H^{\sigma}_{rad}(\Theta)}<\varepsilon\,.
$$
This completes the proof of Lemma~\ref{nedelia}.
\end{proof}
It suffices to prove that
\begin{equation}\label{reduction}
\rho(\Phi(\tau)(K))\geq\rho(K),\quad \tau\in [0,c(1+R)^{-\gamma}],
\end{equation}
where $c$ and $\gamma$ are fixed by Lemma~\ref{nedelia}.
Indeed, it suffices to cover $[0,t]$ by intervals of size $ c(1+R)^{-\gamma}$ and apply
(\ref{reduction}) at each step. Such an iteration is possible since at each
step the image remains a compact of $H^{\sigma}_{rad}(\Theta)$ included in the ball $B_{R}$. 
Let us now prove (\ref{reduction}).
Let $B_{\varepsilon}$ be the open ball in $H^{\sigma}_{rad}(\Theta)$ centered at
the origin and of radius $\varepsilon$.
By the continuity property of $\Phi(t)$, we have that $\Phi(\tau)(K)$ 
is a closed set of $H^{\sigma}_{rad}(\Theta)$
contained in $\Sigma$.
Therefore, by (\ref{close}), we can write
$$
\rho\Big(\Phi(\tau)(K)+\overline{B_{2\varepsilon}}\Big) \geq \limsup_{N\rightarrow \infty}
\rho_{N}\Big(\big(\Phi(\tau)(K)+\overline{B_{2\varepsilon}}\big)\cap E_{N}\Big)\,,
$$
where $\overline{B_{2\varepsilon}}$ is the closed ball in
$H^{\sigma}_{rad}(\Theta)$, centered at the origin and of radius $2\varepsilon$.
Using Lemma~\ref{nedelia}, we obtain that for every $\varepsilon>0$, if we take $N$ large enough, we
have
$$
\big(\Phi_{N}(\tau)(S_{N}(K))+B_{\varepsilon}\big)\cap E_{N}
\subset
\big(\Phi(\tau)(K)+\overline{B_{2\varepsilon}}\big)\cap E_{N}
$$
and therefore
$$
\limsup_{N\rightarrow \infty}
\rho_{N}\Big(\big(\Phi(\tau)(K)+\overline{B_{2\varepsilon}}\big)\cap E_{N}\Big)
\geq 
\limsup_{N\rightarrow \infty}
\rho_{N}\Big(
\big(\Phi_{N}(\tau)(S_{N}(K))+B_{\varepsilon}\big)\cap E_{N}
\Big).
$$
Next, using the uniform continuity property of the flow $\Phi_N$ (see Proposition~\ref{lwpbis}), we obtain that
there exists $c\in]0,1[$, independent of $\varepsilon$ such that for $N$ large enough, we have
$$
\Phi_{N}(\tau)\big((K+B_{c\varepsilon})\cap E_{N}\big)
\subset
\big(\Phi_{N}(\tau)(S_{N}(K))+B_{\varepsilon}\big)\cap E_{N},
$$
where $B_{c\varepsilon}$ is the open ball in $H^{\sigma}_{rad}(\Theta)$ centered at the origin and 
of radius $c\varepsilon$. Therefore
$$
\limsup_{N\rightarrow \infty}
\rho_{N}\Big(\big(\Phi_{N}(\tau)(S_{N}(K))+B_{\varepsilon}\big)\cap E_{N}\Big) 
\geq 
\limsup_{N\rightarrow \infty}\rho_{N}\Big(\Phi_{N}(\tau)\big(
(K+B_{c\varepsilon})\cap E_{N}\big)\Big),
$$
Further, using Proposition~\ref{liouville}, we obtain 
$$
\rho_{N}\Big(\Phi_{N}(\tau)\big((K+B_{c\varepsilon})\cap E_{N}\big)\Big)
=
\rho_{N}\Big((K+B_{c\varepsilon})\cap E_{N}\Big)
$$
and thus
$$
\limsup_{N\rightarrow \infty}\rho_{N}\Big(\Phi_{N}(\tau)\big((K+B_{c\varepsilon})\cap E_{N}\big)\Big)
\geq
\liminf_{N\rightarrow \infty}\rho_{N}\Big((K+B_{c\varepsilon})\cap E_{N}\Big).
$$
Finally, using (\ref{open}), we can write
$$
\liminf_{N\rightarrow \infty}\rho_{N}\Big((K+B_{c\varepsilon})\cap E_{N}\Big)
\geq
\rho(K+B_{c\varepsilon})\geq \rho(K).
$$
Therefore, we have the inequality
$
\rho\Big(\Phi(\tau)(K)+\overline{B_{2\varepsilon}}\Big) \geq \rho(K).
$
By letting $\varepsilon\rightarrow 0$, the dominated convergence gives $\rho(\Phi(\tau)(K))\geq \rho(K)$.
This completes the proof of Theorem~\ref{thm2}.
\end{proof}

\end{document}